\documentclass[1p]{elsarticle}
\usepackage{amssymb}
\usepackage{amsmath}
\usepackage{amsthm}

\usepackage[usenames,dvipsnames]{pstricks}
 \usepackage{epsfig}

\newtheorem{thm}{Theorem}
\newtheorem{lem}[thm]{Lemma}
\newenvironment{pf}{\noindent {\sc Proof.}}{$\Box$\par\medskip}%

\newtheorem{corollary}[thm]{\bf Corollary}%
\newtheorem{example}[thm]{\bf Example}%
\newtheorem{problem}[thm]{\bf Problem}%
\newtheorem{definition}[thm]{\bf Definition}{}%
\newcommand{\be}{\begin{equation}}
\newcommand{\ee}{\end{equation}}
\newcommand{\ol}{\overline}

\newcommand{\wh}{\widehat}

\begin{document}

\begin{frontmatter}

\title{Contractors for flows}

\author{Delia Garijo\fnref{fn1}}
\ead{dgarijo@us.es}
\address{Department of Applied Mathematics I, University of Seville, Seville, Spain}

\author{Andrew Goodall\corref{cor1}\fnref{fn2}} 
\ead{goodall.aj@gmail.com}
\author{Jaroslav Ne\v{s}et\v{r}il\fnref{fn2}}
\ead{nesetril@kam.mff.cuni.cz}
\address{Department of Applied Mathematics and Institute of Theoretical
  Computer Science (ITI), Charles University, Prague, Czech Republic} 
 
\fntext[fn1]{Research supported by projects O.R.I MTM2008-05866-C03-01, PAI FQM-0164 and PAI P06-FQM-01649}

\cortext[cor1]{Corresponding author.}
\fntext[fn2]{Research supported by ITI 1M0545, and the Centre for Discrete
Mathematics, Theoretical Computer Science and Applications (DIMATIA).}




\begin{abstract}
We answer a question raised by Lov\'{a}sz and B.~Szegedy [Contractors and connectors in graph algebras, {\em J. Graph Theory} 60:1 (2009)] asking for a contractor for the graph parameter counting the number of $B$-flows of a graph, where $B$ is a subset of a finite Abelian group closed under inverses. We prove our main result using the duality between flows and tensions and finite Fourier analysis. 
We exhibit several examples of contractors for $B$-flows, which are of interest in relation to the family of $B$-flow conjectures formulated by Tutte, Fulkerson, Jaeger, and others.
\end{abstract}
\end{frontmatter}

\section{Introduction}

In their paper~\cite{LS09} Lov\'{a}sz and B.~Szegedy introduce the notion of contractors and connectors in quantum graph algebras. Connectors and contractors for a graph parameter allow expressions for the parameter of a graph obtained by the connection of two vertices by an edge or the contraction of an edge (either operation possibly producing multiple edges or loops, respectively) in terms of values of the parameter on graphs obtained by operations that preserve simplicity. For example, a connector for the flow polynomial is any path of length two or more.
Contractors for graph parameters can be thought of as a generalization of the deletion-contraction identity for the chromatic polynomial, flow polynomial, and other specializations of the Tutte polynomial. 

In~\cite[Section 2.3]{LS09} Lov\'asz and B.~Szegedy comment that ``there does not seem to be a simple explicit construction for a contractor'' for the number of $B$-flows, where $B$ is a subset of a finite Abelian group closed under inverses. In Section~\ref{sec:contractor_Sflows} we give an explicit construction (Theorem~\ref{thm.S-flows} and its Corollary~\ref{cor:disconn}), followed by numerous examples of relevance to the long-standing conjectures (such as the Cycle Double Cover Conjecture and Fulkerson's Conjecture) that Jaeger reformulates as a unifying $B$-flow conjecture in~\cite{Jae88}. 

We prove our result using the Fourier transform on finite Abelian groups. Towards this end, in
Section~\ref{sec:prelim} we introduce homomorphism functions, flows and tensions of a graph, connectors and contractors in graph algebras, and then describe those properties of the Fourier transform that will be used to enable us to reach the main result (Theorem~\ref{thm:contractor from Fourier transform}, and its special case for $B$-flows, Theorem~\ref{thm.S-flows}).
The Fourier transform expresses any weighted function on the set of tensions of a graph (effectively a homomorphism function to an edge-weighted graph) as a related weighted function on the set of flows of the graph. This in particular allows the number of $B$-flows to be written as the  number of  homomorphisms to an edge-weighted graph (Lemma~\ref{lem:Sflows}). Although the duality relationship between flows and tensions is well known, there have not been many useful applications of this correspondence. In this paper it allows us to transform a special type of instance of a problem we do not in general know yet how to solve (Problem~\ref{prob:contractor_eigenvalue_zero} in Section~\ref{sec:conclusion}) into one that we can solve.    

For graphs $G$ and $H$, the number of homomorphisms from $G$ to $H$ is denoted by ${\rm hom}(G, H)$, and we extend this notation to edge-weighted graphs $H$. (A precise definition is given in Section~\ref{sec:hom}). In Section~\ref{sec:connectors_contractors} we consider the problem of finding a connnector and contractor for the homomorphism number ${\rm hom}(\,\cdot\,, H)$.
The minimum polynomial of the adjacency matrix of $H$ can be used to construct a connector for ${\rm hom}(\,\cdot\,, H)$ , and also a contractor for ${\rm hom}(\,\cdot\,, H)$ if the adjacency matrix of $H$ does not have eigenvalue $0$ (Theorem~\ref{thm:connector contractor}).
The relevant edge-weighted graph $H$ that gives the number of $B$-flows turns out to have eigenvalue $0$, but its special structure as an edge-weighted Cayley graph allows us to obtain a contractor by using the minimum polynomial of the graph whose edge weights are obtained by taking the Fourier transform of the edge weights of $H$ (Theorem~\ref{thm:contractor from Fourier transform}).
In Section~\ref{sec:contractor_Sflows} we apply the theory developed in Sections~\ref{sec:prelim} and~\ref{sec:connectors_contractors} to derive our main result, giving a contractor for the number of $B$-flows (Theorem~\ref{thm.S-flows}).

As an application of this result, in Section~\ref{sec:ex} we give a series of examples chosen on account of their relevance to such problems as the Cycle Double Cover Conjecture, Fulkerson's Conjecture and the Petersen Flow Conjecture, where we follow Jaeger's formulation~\cite{Jae88} of these conjectures as $B$-flow problems. We also consider the example $B=\{-1,+1\}$ in the additive group $\mathbb{Z}_n$ of integers modulo $n$, because the number of $B$-flows here is equal to the number of Eulerian orientations of $G$ when the vertex degrees of $G$ all belong to $\{0,1,\ldots, n\!-\!1,n\!+\!1\}$. As shown in~\cite{LS06}, the number of Eulerian orientations of a graph is known not to be a homomorphism number to a finite edge-weighted graph, and therefore does not have a contractor~\cite[Example 2.6]{LS09}. The interesting question arises of whether there is a ``limiting contractor'' for Eulerian orientations.

In Section~\ref{sec:conclusion} we conclude by posing some further problems concerning contractors of homomorphism functions to finite edge-weighted graphs.

\section{Preliminaries}\label{sec:prelim}

\subsection{Homomorphism functions}\label{sec:hom}

Given two graphs $G=(V(G),E(G))$ and
$H=(V(H),E(H))$, a \emph{homomorphism} of $G$ to  $H$,
written as $\psi:G\rightarrow H$, is a mapping $\psi:V(G)\rightarrow
V(H)$ such that $\psi(u)\psi(v)\in E(H)$ whenever $uv\in E(G)$. The number of homomorphisms of $G$ to $H$ is denoted by hom$(G,H)$. For a fixed graph $H$, hom$(\,\cdot\,, H)$ is a \emph{graph parameter}, i.e., a function defined on graphs invariant under isomorphism. This parameter can be generalized for weighted graphs $H$ as follows. 

Let $H$ be a weighted graph with a real positive weight $\alpha_H(i)$ associated with each vertex $i\in V(H)$ and a real weight $\beta_H(ij)$ associated with each edge $ij\in E(H)$. Let $G$ be an unweighted graph and fix
$U\subseteq V(G)$. Define for every function
$\phi:U\rightarrow V(H)$ the weight 
$${\rm hom}_\phi(G,H)=\sum_{\stackrel{\psi:V(G)\rightarrow
    V(H)}{\psi_{|U}=\phi}}\prod_{u\in V(G)\setminus U}\alpha_H(
\psi(u)) \prod_{uv\in E(G)}\beta_H
(\psi(u),\psi(v)).$$

Then the number of  homomorphisms of $G$ to $H$ is defined by
 
$${\rm hom}(G,H)=\sum_{\phi:U\rightarrow V(H)}\left(\prod_{u\in  U}\alpha_H(
\phi(u))\right){\rm hom}_\phi(G,H).$$

In this paper we assume $ \alpha_H(i)=1$ for every $i\in V(H)$ (vertex weights equal to one)  as
   it simplifies the
  exposition  and also because for positive integer
  vertex weights we can replace the vertex weight $ \alpha_H(i)$    by splitting the
  vertex $i$ into  $ \alpha_H(i)$         unweighted twin vertices. Thus, $H$ will be considered as  an \emph{edge-weighted graph} with symmetric adjacency matrix $A=A(H)=(\beta_H(ij))$ and

$${\rm hom}_\phi(G,H)=\sum_{\stackrel{\psi:V(G)\rightarrow
    V(H)}{\psi_{|U}=\phi}}  \prod_{uv\in E(G)}\beta_H
(\psi(u)\psi(v)),$$

$${\rm hom}(G,H)=\sum_{\phi:U\rightarrow V(H)}{\rm hom}_\phi(G,H).$$

\subsection{Flows and tensions}\label{sec:flows_tensions}

Let $G=(V,E)$ be a graph with an arbitrary fixed orientation of its edges. Denote by $c(G)$ the number of connected components of $G$, and by $r(G)=|V|-c(G)$ its rank. Let $D$ be the $V\times E$ incidence matrix of $G$, with $(v,e)$-entry given by
$$D_{v,e}=\begin{cases} +1 & \mbox{\rm $v$ is the head of $e$,}\\
-1 & \mbox{\rm $v$ is the tail of $e$,}\\
0 & \mbox{\rm $v$ is not an endpoint of $e$, or $e$ is a loop on $v$.}\end{cases}$$

Suppose $\Gamma$ is a finite commutative ring with unity, and let $\Gamma^V=\{f:V\rightarrow\Gamma\}$ and $\Gamma^E=\{g:E\rightarrow \Gamma\}$. 

The linear transformation $D:\Gamma^E\rightarrow\Gamma^V$ of $\Gamma$-modules has kernel equal to the set of {\em $\Gamma$-flows} of $G$. This is to say a function $g:E\rightarrow \Gamma$ is a $\Gamma$-flow if and only if, for each vertex $v\in V$, 
$$\sum_{\stackrel{e=uv}{\mbox{\tiny $v$ head of $e$}}}g(e)-\sum_{\stackrel{e=uv}{\mbox{\tiny $v$ tail of $e$}}}g(e)=0.$$ 

If $B=-B$ is a subset of $\Gamma$ then an {\em $B$-flow} of $G$ is a $\Gamma$-flow taking values in $B$. (We have followed Jaeger's~\cite{Jae88} usage of the term $B$-flow; in~\cite{LS09} Lov\'asz and B.~Szegedy use the letter $S$ instead of $B$, and so speak of $S$-flows rather than the $B$-flows.)
In particular, when $B=\Gamma\setminus\{0\}$ an $B$-flow is a {\em nowhere-zero $\Gamma$-flow} of $G$. It is well known that the number of nowhere-zero $\Gamma$-flows of $G$ depends only on $|\Gamma|$ and is given by the evaluation of the {\em flow polynomial} $F(G;z)$ at $z=|\Gamma|$.

 The transpose linear transformation $D^\top:\Gamma^V\rightarrow \Gamma^E$ 
has image the set of {\em $\Gamma$-tensions} of $G$. A $\Gamma$-tension $g:E\rightarrow \Gamma$ arises from a vertex $\Gamma$-colouring $f:V\rightarrow \Gamma$ of $G$, by setting 
$g(e)=f(v)-f(u)$ for each edge $e$ with tail $u$ and head $v$.

A {\em $B$-tension} of $G$ is a $\Gamma$-tension of $G$ taking values in $B$. In particular, when $B=\Gamma\setminus\{0\}$ a $B$-tension arises from $|\Gamma|^{c(G)}$ proper vertex $\Gamma$-colourings of $G$ (a proper colouring by definition has the property that the difference between endpoint colours on any edge is non-zero). The number of proper vertex $\Gamma$-colourings is given by the evaluation of the {\em chromatic polynomial} $P(G;z)$ at $z=|\Gamma|$.

Note that to define $\Gamma$-flows and $\Gamma$-tensions all that is needed is the structure of $\Gamma$ as an additive Abelian group. As we shall see in Section~\ref{sec:fourier}, it is convenient to impose the extra ring structure on $\Gamma$ in order to relate the sets of tensions and flows to each other as orthogonal complements. 

\subsection{Connectors and contractors}

For fixed positive integer $k$,
a {\em $k$-labelled graph} $(G,\lambda)$ is a finite graph $G$ (without loops but possibly multiple edges)              
 together with a
function $\lambda:[k]\rightarrow V(G)$. Lov\'asz and
Szegedy~\cite{LS09} restrict $\lambda$ to be an injective function,
but we shall follow Schrijver~\cite{Schr09} and allow vertices to have
multiple labels.
A $k$-labelled graph $(G,\lambda)$ is {\em simple} if $G$ has no
parallel edges, the function $\lambda$ is
injective and its image is an independent set of vertices. 

A $k$-labelled {\em quantum graph} is a formal linear combination of
$k$-labelled graphs with coefficients in $\mathbb{R}$. The set of
$k$-labelled quantum graphs is denoted by $\mathcal{G}_k$. The subset
of $k$-labelled quantum graphs whose labelled vertices are independent
is denoted by $\mathcal{G}_k^0$, and the subset of simple $k$-labelled
quantum graphs is denoted by $\mathcal{G}^{\mbox{\rm \tiny simp}}_k$.


We let $K_n$ denote the complete graph on $n$ vertices and $\ol{K}_n$ its complement, consisting of $n$ isolated vertices. 

For two $k$-labelled graphs $X$ and $Y$, the {\em product} $XY$ is
defined~\cite{LS09} by taking the disjoint union of $X$ and $Y$ and then
identifying vertices which share the same label. The product is
associative and commutative, and extends linearly to
$\mathcal{G}_k$. The identity for this multiplication on
$\mathcal{G}_k$ is $(\ol{K}_k,\lambda)$ where $\lambda$ is injective,
i.e., $k$ independent vertices each with one label.

In the sequel we shall just consider the case $k=2$.

A further binary operation on $2$-labelled graphs is the {\em concatenation} $X\circ Y$ of $X$ and $Y$, defined by identifying the
vertex of $X$ labelled $2$ with the vertex labelled $1$ in $Y$ and
unlabelling this merged vertex. Concatentation is associative but not
commutative. It is extended linearly to $\mathcal{G}_2$. The identity
for this multiplication on $\mathcal{G}_2$ is $K_1$, with
both labels $1$ and $2$ on the single vertex.

Let $P_\ell$, $\ell\geq 1$, denote the $2$-labelled path on $\ell+1$ vertices with one endpoint
labelled $1$ and the other endpoint labelled $2$. In this notation, $P_1=K_2$. We have $P_\ell\circ
P_m=P_{\ell+m}$, and $P_\ell P_m=C_{\ell+m}$, where the cycle has
labels spaced $\ell$ edges apart in one direction, $m$ edges apart
in the other. 
If $X$ is a $2$-labelled graph then $XP_1$ is the
$2$-labelled graph with an edge added joining the vertices labelled
$1$ and $2$. In particular, if $X$ is a simple $2$-labelled graph then
$XP_1$ joins non-adjacent vertices with an edge.

The graph $K_1$ with vertex carrying both labels $1$ and $2$ when
applied as a product with $X$ (giving the graph $XK_1$) identifies the vertices labelled $1$ and
$2$ in $X$ (or does nothing if both labels are on a single vertex in $X$).  
If $1$ and $2$ are independent in $X$, then $XK_1$ is the graph $XP_1$
with the edge joining $1$ and $2$ contracted.

 The graph $\ol{K}_2$ with one vertex labelled 1 and the other labelled 2  is the identity for the product of 2-labelled
graphs. 

A {\em series-parallel} $2$-labelled graph is obtained from $P_1$ by a
sequence of products (parallel extensions) and concatenations (series
extensions). 

\begin{definition}{\rm \cite{LS09}}
Let $h$ be a graph parameter. A $2$-labelled quantum graph $Z\in\mathcal{G}_2$ is a
{\em contractor} for $h$ if for all $X\in \mathcal{G}_2^0$ we have
$$h(XZ)=h(XK_1).$$

A simple $2$-labelled quantum graph $Z\in \mathcal{G}^{\mbox{\rm \tiny
    simp}}_k$ is a {\em connector} for $h$ if for all
    $X\in\mathcal{G}_2$ we have
$$h(XZ)=h(XP_1).$$

\end{definition}

Note that attaching a contractor $Z$ at two non-adjacent vertices acts like identifying these vertices, and attaching a connector $Z$ at two vertices acts like joining them with an edge; if
these vertices are already joined by an edge then since
$Z\in\mathcal{G}_2^{\mbox{\rm \tiny simp}}$ the graph $XZ$ remains
simple whereas $XP_1$ is not simple.

\begin{example}
The deletion-contraction relation for the chromatic polynomial takes the
form $P(XP_1;z)=P(X;z)-P(XK_1;z)$ where $X$ has its two
labels on non-adjacent vertices. Extending the chromatic polynomial linearly
to $\mathcal{G}_2$, this is to say $P(XK_1;z)=P(X(\ol{K}_2-P_1);z)$. The
quantum graph $\ol{K}_2-P_1$ is a contractor for the chromatic
polynomial. 

For $X\in\mathcal{G}_2^{\mbox{\rm \tiny simp}}$, the flow polynomial
satisfies $F(XP_1;q)=F(XP_\ell;q)$ for any $\ell\geq 2$; the
2-labelled graph $P_2$
is a simple connector for the flow polynomial. 
\end{example}


\subsection{The Fourier transform}\label{sec:fourier}

The duality between flows and tensions is well known, especially in statistical physics, see e.g.~\cite{B77}, with an early example being Van der Waerden's Eulerian subgraph expansion of the Ising model partition function~\cite{vdW41}. 
The material in this section is standard and can be found in~\cite{Terras99}, except for its main result, Lemma~\ref{lem:Sflows}, for which we therefore provide a proof. We briefly recall the relevant facts here for the convenience of the reader. More examples of applications of elementary Fourier analysis on finite Abelian groups to graph theory can be found in~\cite{CCC07}.

Let $\Gamma$ be a finite additive Abelian group of order $n$.
If $\ell$ is the least common multiple of the orders of the elements of $\Gamma$ then there is a unique sequence of positive integers $n_1,n_2\ldots, n_r$ such that $2\leq n_1\;|\;n_2\;|\cdots |\;n_r=\ell$ and 
\be\label{eq:decomp}\Gamma\cong\mathbb{Z}_{n_1}\oplus\mathbb{Z}_{n_2}\oplus\cdots \oplus\mathbb{Z}_{n_r}.\ee
We can give $\Gamma$ the structure of a ring by defining multiplication of elements $x=(x_1.x_2,\ldots, x_r)$ and $y=(y_1,y_2,\ldots, y_r)$ componentwise: $xy=(x_1y_1,x_2y_2,\ldots, x_ry_r)$. 
Alternative multiplicative structures are possible. For example, when $\Gamma$ is isomorphic to $\mathbb{Z}_p\oplus \mathbb{Z}_p\oplus\cdots\oplus\mathbb{Z}_p$ we can endow $\Gamma$ with the structure of the finite field $\mathbb{F}_{p^r}$. 
A {\em character} of $\Gamma$ is a homomorphism $\chi:\Gamma\rightarrow\mathbb{C}^{\times}$, where $\mathbb{C}^\times$ is the multiplicative group of the field of complex numbers. The set of characters of $\Gamma$ form a group $\wh{\Gamma}$ under pointwise multiplication, which is isomorphic to $\Gamma$ when $\Gamma$ is a finite Abelian group.
For each $x\in\Gamma$, let $\chi_x$ denote the image of $x$ under a fixed isomorphism $\Gamma\rightarrow\wh{\Gamma}$. In particular, the trivial character $\chi_0$ is defined by $\chi_0(y)=1$ for all $y\in \Gamma$, and $\chi_{-x}(y)=\ol{\chi_x(y)}$, the bar denoting complex conjugation.

Suppose then that the additive group $\Gamma$ has been given the extra structure of a commutative ring. 
A character $\chi\in\wh{\Gamma}$ is a {\em generating character} for $\Gamma$ if $\chi_x(y)=\chi(xy)$ for each character $\chi_x\in\wh{\Gamma}$.
For example, the ring $\mathbb{Z}_n$ has a generating character defined by $\chi(x)=e^{2\pi i x/n}$, and the field $\mathbb{F}_{p^r}$ has a generating character defined by $\chi(x)=e^{2\pi\mbox{\rm \tiny Tr}(x)/p}$, where ${\rm Tr}(x)=x+x^p+\cdots + x^{p^{r-1}}$. 
For any additive Abelian group $\Gamma$ we can use its canonical expression~\eqref{eq:decomp} as a direct sum of rings of the form $\mathbb{Z}_n$ to obtain a generating character. 

Consider $\Gamma^m$, the $m$-fold direct sum of $\Gamma$, and suppose $\psi$ is a
generating character for $\Gamma$. Then $\chi$ defined by
$\chi(x_1,\ldots,x_m)=\psi(x_1)\cdots\psi(x_m)$ for $(x_1,\ldots,
x_m)\in \Gamma^m$ is a generating character for $\Gamma^m$. The {\em Euclidean
inner product} (dot product) is defined for $x=$$(x_1,\ldots, x_m),
y$$=(y_1,\ldots, y_m)\in \Gamma^m$ by $x\cdot y$$=x_1y_1+\cdots +x_my_m$. Since
$\psi(x_1)\cdots\psi(x_m)=\psi(x_1+\cdots x_m)$, it follows that $\chi_x(y)=\chi(xy)=\psi(x\cdot y)$ for $x,y\in
\Gamma^m$.

Denote by
$\mathbb{C}^\Gamma$ the vector space over $\mathbb{C}$ of all functions from $\Gamma$ to $\mathbb{C}$. This is an inner product space with Hermitian
inner product $\langle\,\cdot\,\rangle$ defined for
$\alpha,\beta\in\mathbb{C}^\Gamma$ by
$$\langle \alpha,\beta\rangle=\sum_{x\in \Gamma}\alpha(x)\ol{\beta(x)}.$$


The vector space $\mathbb{C}^\Gamma$ has the additional structure of an
algebra under either of the following two definitions of
multiplication:

\begin{itemize}
\item[(i)]{ the {\em pointwise product} $\alpha\cdot \beta$ of
    $\alpha,\beta\in\mathbb{C}^\Gamma$, defined for $x\in \Gamma$ by
$\alpha\cdot \beta(x)=\alpha(x)\beta(x)$,}
\item[(ii)] {the {\em convolution} $\alpha\ast \beta$ of $\alpha,\beta\in\mathbb{C}^\Gamma$,
    defined for $x\in \Gamma$ by
$$\alpha\ast \beta(x)=\sum_{y\in \Gamma}\alpha(y)\beta(x-y).$$}
\end{itemize}


The set $\{\delta_x:x\in \Gamma\}$ of indicator functions defined by 
$$\delta_x(y)=\begin{cases}
1 & x=y,\\
0 & x\neq y,\end{cases}$$
form an orthonormal basis for $\mathbb{C}^\Gamma$, with $\langle \delta_x,\delta_y\rangle=\delta_x(y)$. We extend the indicator function
notation to subsets $B$ of $\Gamma$, defining $\delta_B=\sum_{x\in B}\delta_x$.

The characters of $\Gamma$ are also orthogonal in this inner product space, with
$\langle\chi_x,\chi_y\rangle=|\Gamma|\delta_x(y)$. 

Fix an isomorphism $x\mapsto \chi_x$ of $\Gamma$ with $\wh{\Gamma}$ and let
$\chi$ be a generating character for $\Gamma$ such that
$\chi_x(y)=\chi(xy)$.

For $\alpha\in\mathbb{C}^{\Gamma}$ the {\em Fourier transform}
$\wh{\alpha}\in\mathbb{C}^{\Gamma}$ is defined for $y\in \Gamma$ by 
\be\label{FT}{\wh{\alpha}(y)=\langle \alpha,\chi_y\rangle=\sum_{x\in
  \Gamma}\alpha(x)\chi_y(-x).}\ee

The Fourier transform maps the basis of indicator functions to the
basis of characters: $\wh{\delta_y}=\chi_{-y}$.
The {\em Fourier inversion
  formula} $\wh{\wh{\alpha}}(x)=|\Gamma|\alpha(-x)$,
gives the inverse transform
\be\label{inversion}{\alpha(x)=\frac{1}{|\Gamma|}\langle \wh{\alpha},\chi_{-x}\rangle=\frac{1}{|\Gamma|}\sum_{y\in \Gamma}\wh{\alpha}(y)\chi_x(y).}\ee
Note that $\wh{\delta_\Gamma}=|\Gamma|\delta_0$ and $\wh{\delta_0}=\delta_\Gamma$, since $\langle \delta_\Gamma,\chi_y\rangle=\langle\chi_0,\chi_y\rangle=|\Gamma|
\delta_0(y)$. 




The Fourier transform gives an isomorphism of the
algebra $\mathbb{C}^\Gamma$ with multiplication pointwise product with the algebra
$\mathbb{C}^\Gamma$ with multiplication convolution: for $y\in \Gamma$,
\be\label{convolution}{\wh{\alpha\ast \beta}(y)=\wh{\alpha}\cdot\wh{\beta}(y),}\ee
\be\label{pointwise}{\wh{\alpha\cdot
    \beta}(y)=\frac{1}{|\Gamma|}\wh{\alpha}\ast\wh{\beta}(y).}\ee

For a subgroup $C$ of $\Gamma$, the {\em annihilator}
$C^{\sharp}$ of
$C$ is defined by
$$C^{\sharp}:=\{x\in R:\chi_x(y)=1\;\mbox{\rm for all }y\in C\}.$$
The annihilator $C^\sharp$ is a subgroup of $\Gamma$ isomorphic to $\Gamma/C$.
A second key property of the Fourier transform is that it takes indicators of subgroups to (scalar
multiples of) indicators of
their annihilators:
\be\label{FTannihil}{\wh{\delta}_C(y)=\sum_{x\in C}\chi_x(y)=|C|1_{C^\sharp}(y).}\ee

By (\ref{inversion}), (\ref{convolution}) and
(\ref{FTannihil}) there follows the {\em Poisson summation formula},
\be\label{Poisson}{\sum_{x\in C}\alpha(x+z)=\frac{1}{|C^\sharp|}\sum_{x\in
    C^\sharp}\wh{\alpha}(x)\chi_{z}(x).}\ee


Consider again the $m$-fold direct sum  $\Gamma^m$.
The {\em orthogonal} to a subset $C$ of $\Gamma^m$ (with respect
to the Euclidean inner product) is defined by
\[C^\perp=\{y\in \Gamma: x\cdot y=0\;\mbox{\rm for all }y\in C\}.\]
If $\Gamma$ has a
    generating character and $C$ is a $\Gamma$-submodule of $\Gamma^m$,  
then $C^\sharp=C^\perp$.

Thus we have the following special case of the Poisson summation formula when $C$ is a $\Gamma$-submodule of $\Gamma^m$ with orthogonal complement $C^\perp$:
\be\label{eq:poisson} \sum_{x\in C}\alpha(x)=\frac{1}{|C^\perp|}\sum_{y\in C^\perp}\wh{\alpha}(y).\ee





We now apply the machinery we have set up to relate flows and tensions of a graph.

Let $D$ be the incidence matrix of a graph $G$, as defined in Section~\ref{sec:flows_tensions}. Recall that the kernel $\ker D$ of the map $D:\Gamma^E\rightarrow \Gamma^V$ consists of the $\Gamma$-flows of $G$, and the image ${\rm im}D^\top$ of the adjoint map $D^\top:\Gamma^V\rightarrow \Gamma^E$ consists of the $\Gamma$-tensions of $G$.



Provided $\Gamma$ has a generating character, we have
  $({\rm im}\, D^\top)^\perp=\ker D$ and $(\ker D)^{\perp}={\rm im\,}D^\top$.
We also have $|\ker D|=|\Gamma|^{|E|-r(G)}$ and $|{\rm im}\, D^\top|=|\Gamma|^{r(G)}$.
Hence by~\eqref{eq:poisson},
\be\label{eq:flow_tensions}\sum_{x\in\ker D}\alpha(x)=\frac{1}{|\Gamma|^{r(G)}}\sum_{y\in\mbox{\rm \tiny im}D^\top}\wh{\alpha}(y).\ee

\begin{lem}\label{lem:Sflows}
Let $B=-B$ be a subset of an additive Abelian group $\Gamma$. Let $H={\rm Cayley}(\Gamma, B)$ be the graph on vertex set $\Gamma$ with an edge joining vertices $i$ and $j$ if and only if $j-i\in B$, and let $\wh{H}$ be the edge-weighted Cayley graph on vertex set $\Gamma$ with edge $ij$ having weight $\wh{\delta}_B(j-i)$.
 
Then the number of $B$-tensions of a graph $G=(V,E)$ is equal to\newline $|\Gamma|^{-c(G)}{\rm hom}(G,H)$ and the number of $B$-flows of $G$ is equal to 
$|\Gamma|^{-|V|}{\rm hom}(G,\wh{H})$.  
\end{lem}

\begin{pf}
A function $f:V\rightarrow\Gamma$ is a homomorphism contributing to ${\rm hom}(G,H)$ if and only if $f(v)-f(u)\in B$ for each edge $e=uv\in E$, which is to say that $g(e):=f(v)-f(u)$ is a $B$-tension of $G$. For a given $B$-tension $g:E\rightarrow \Gamma$ there are $|\Gamma|^{c(G)}$ vertex colourings $f:V\rightarrow\Gamma$ with $Df=g$.
This proves that the number of $B$-tensions is $|\Gamma|^{-c(G)}{\rm hom}(G,H)$. 
 
To prove that $|\Gamma|^{-|V|}{\rm hom}(G,\wh{H})$ is the number of $B$-flows, take $\alpha=\delta_{B^E}$ in~\eqref{eq:flow_tensions}, which has Fourier transform for $g:E\rightarrow\Gamma$ given by
$$\wh{\delta}_{B^E}(g)=\prod_{e\in E}\wh{\delta}_B(g(e)),$$
which yields
$$|\ker D\cap B^E|=\frac{1}{|\Gamma|^{|V|}}\sum_{f:V\rightarrow\Gamma}\prod_{uv\in E}\wh{\delta_B}(f(v)-f(u)).$$
\end{pf}

\section{Connectors and contractors for homomorphism numbers}\label{sec:connectors_contractors}

After the preliminaries of Section~\ref{sec:prelim} we now turn to the general problem of finding a contractor for homomorphism numbers ${\rm hom}(\,\cdot\,, H)$. (In Section~\ref{sec:contractor_Sflows} we then consider the special case of finding a contractor for the number of $B$-flows.) In this section we shall also find an expression for a connector for ${\rm hom}(\,\cdot\,, H)$, although this is of less interest for the case of $B$-flows (where $P_2$ is always a connector).

Let $H$ be an edge-weighted graph with adjacency matrix $A=(\beta(ij))$. To every $X\in\mathcal{G}_2$ assign a $V(H)\times V(H)$-matrix $M(X)=M_H(X)$
with $(i,j)$-entry equal to ${\rm hom}_\phi(X,H)$ where $\phi(1)=i$
and $\phi(2)=j$. As $H$ is fixed we shall usually leave the dependence on $H$ implicit and write $M(X)$ instead of $M_H(X)$.  The function
$M:\mathcal{G}_2\rightarrow\mathbb{C}^{V(H)\times V(H)}$ is linear and
satisfies
$$M(X\circ Y)=M(X)M(Y), \;\;M(XY)=M(X)\circ M(Y),$$
where $\circ$ is the entrywise (Hadamard, or Schur) product.

Note that $M(P_1)=(\beta(ij))$ is the adjacency matrix of $H$. 
Let $C_k^*$ denote the graph on two vertices joined by $k$ parallel edges (the dual of the $k$-cycle), with one vertex labelled $1$ and the other vertex labelled $2$. The graph $P_k^*$ consists of one vertex with $k$ loops attached (the dual of the path on $k$ edges), and with labels $1$ and $2$ on the single vertex. 
\begin{lem}\label{lem:path and thick edge}
If $H$ has adjacency matrix $A$ then
$$M(P_k)=A^{k}, \;\;\; M(C_k^*)=A^{\circ k},$$
where $A^{\circ k}$ denotes the $k$-fold Hadamard product $A\circ
A\circ\cdots\circ A$, and $A^k$ the $k$-fold ordinary matrix product. 
\end{lem}

\begin{pf} Since $M(P_1)=A$, $P_k=P_1^{\circ k}$ and $C_k^*=P_1^k$ then $M(P_k)=M(P_1^{\circ k})=M(P_1)^k=A^k$ and 
 $M(C_k^*)=M(P_1^{ k})=M(P_1)^{\circ k}=A^{\circ k}$.\end{pf}

Lemma~\ref{lem:path and thick edge} includes for $k=0$ the fact that if $K_1$ is single vertex with
labels $1$ and $2$ then $M(K_1)=I$. 
We also have $M(\ol{K}_2)=J$, where $\ol{K}_2$ comprises two isolated vertices, one labelled $1$
the other $2$. (The weight of a homomorphism from a graph $X$ with no edges to $H$ is equal to the number of maps from $V(X)$ to $V(H)$.
Restricting to homomorphisms coinciding with $\phi$ on the labelled vertices $1$ and $2$ of $X$, we have ${\rm hom}_\phi(X,H)=1$ when the vertices labelled $1$ and $2$ are distinct, and ${\rm hom}_\phi(X,H)=\delta_{\phi(1),\phi(2)}$ when the same vertex of $X$ has labels $1$ and $2$.)

A generalization of Lemma~\ref{lem:path and thick edge} is that the matrices generated by $A$ (by applying 
the ordinary and Hadamard products of matrices) are equal to $M(X)$ for 
 a series-parallel graph $X$.

\begin{lem}
Let $C_k$ denote the cycle on $k$ vertices with one vertex receiving
both labels $1$ and $2$, and $P_k^*$ the graph comprising $k$ loops on one
vertex with both labels
$1$ and $2$. Then
$$M(C_k)=I\circ A^k,\;\; M(P_k^*)=I\circ A^{\circ k}.$$
\end{lem}
\begin{pf}
The result holds by Lemma \ref{lem:path and thick edge}, since $C_k=P_kK_{1}$ and $P_k^*=C_k^*K_1$.
\end{pf}

A fundamental observation with regard to connectors and contractors for homomorphism numbers is the following:
\begin{lem} {\rm \cite{LS09}}
A 2-labelled quantum graph $Z$ is a contractor for ${\rm hom}(\; .\; ,
H)$ if and only if $M(Z)=I$ and a simple 2-labelled quantum graph $Z$
is a connector for ${\rm hom}(\; .\;, H)$ if and only if $M(Z)=A$.
\end{lem}

Using this lemma we can immediately write down a connector for ${\rm hom}(\,\cdot\,, H)$, and also a contractor when the adjacency matrix of $H$ does not have $0$ as an eigenvalue. (Compare~\cite[Section 4.1]{LS09}.)
\begin{thm}\label{thm:connector contractor}
Suppose the adjacency matrix $A$ of $H$ has minimum polynomial 
$$p_A(t)=p_0+p_1t+\cdots + p_\ell t^\ell.$$
\begin{itemize}
\item[(i)]
If $p_0\neq 0$ (i.e, if $0$ is not an eigenvalue of $A$), then a connector for ${\rm hom}(\,\cdot\,,H)$ is given by
$$Z=-\frac{1}{p_0}\left[p_1P_2+p_2P_3+\cdots + p_\ell P_{\ell+1}\right],$$
where $P_{k}=P_1^{\circ k}$ is the path with $k$ edges, and with one endpoint
vertex labelled $1$, the other labelled $2$. Also, 
a contractor for ${\rm hom}(\,\cdot\,,H)$ is given by
$$Z=-\frac{1}{p_0}\left[p_1P_1+p_2P_2+\cdots +p_{\ell}P_{\ell}\right].$$

\item[(ii)] If $p_0=0$, then $p_1\neq 0$ and 
 a
connector for ${\rm hom}(\,\cdot\, ,H)$ is given by 
$$Z=-\frac{1}{p_1}\left[p_2P_2+\cdots + p_\ell P_{\ell}\right].$$
\end{itemize}
\end{thm}
\begin{pf}
Note that $P_{k}=P_1^{\circ k}$ and $M(P_{k})=A^k$. 
Suppose first $p_0\neq 0$. Then $I=-\frac{1}{p_0}\sum_{k=1}^\ell p_kA^k.$
Hence, if $Z=-\frac{1}{p_0}\sum_{k=1}^\ell p_k P_{k}$ then $M(Z)=I$, and $Z$ is a contractor for ${\rm hom}(\,\cdot\,, H)$.
Likewise,
$A=-\frac{1}{p_0}\sum_{k=1}^\ell p_kA^{k+1}$
so that $Z=-\frac{1}{p_0}\sum_{k=1}^\ell p_kP_{k+1}$ is a connector for ${\rm hom}(\,\cdot\,, H)$.

When $p_0=0$ and $p_1\neq 0$, we have
$A=-\frac{1}{p_1}\sum_{k=2}^\ell p_kA^{k}$ from which it follows that $Z=-\frac{1}{p_1}\sum_{k=2}^\ell p_kP_{k}$ is a connector.
\end{pf}

We are left with the problem of constructing a contractor for ${\rm hom}(\,\cdot\,, H)$ when the adjacency matrix of $H$ has $0$ as an eigenvalue. 
When $H$ is a weighted Cayley graph we are able to find such a contractor, and the remainder of this section is devoted to showing how.

 We begin with the following result of Hoffmann, for which see e.g. \cite[Chapter 3]{B94}.
\begin{thm} Let $A$ be
  the adjacency matrix of a graph $H$. The all-one matrix $J$ is a polynomial
  in $A$ if and only if $H$ is connected and $A$ has eigenvector
  $\mathbf{1}$ (i.e., $H$ is regular). 
\end{thm}
%

\begin{corollary}\label{cor:J in powers of A}
Suppose that the adjacency matrix $A$ of a connected graph $H$ on $n$
vertices has eigenvector $\mathbf{1}$ with eigenvalue $s$, that its
distinct eigenvalues are $s>\lambda_1>\cdots >\lambda_{r-1}$, and that
$$q(t)=\prod_{1\leq i\leq r-1}(t-\lambda_i).$$
Then 
$$J=\frac{n}{q(s)}q(A).$$
\end{corollary}
\begin{pf}
We have $q(A)=\ell J$ for some $\ell\neq 0$ and $q(A)$ has eigenvalues
$q(s)$ and $q(\lambda_i)=0$ for $1\leq i\leq r-1$. Since the only
non-zero eigenvalue of $\ell J$ is $\ell n$ it follows that $\ell n=q(s)$.
\end{pf}
Note that the polynomial $q(t)$ in Corollary~\ref{cor:J in powers of A} is equal to $p_A(t)/(t-s)$, where $p_A(t)$ is the minimum polynomial of $A$. (Since $A$ is symmetric it is diagonalizable and so each of its Jordan blocks has size $1$, which implies the multiplicity of each root of the minimum polynomial of $A$ is equal to one.)  

Suppose now that $A$ has rows and columns indexed by an additive Abelian
group $\Gamma$ of order $n$ and has $(i,j)$ entry equal to $(\beta(i-j))$ for all
$i,j\in\Gamma$, where $\beta:\Gamma\rightarrow\mathbb{C}$ satisfies
$\beta(-i)=\beta(i)$. Thus $A$ is the adjacency matrix of a weighted Cayley graph on $\Gamma$.

Recall the definition of characters and the Fourier transform from Section~\ref{sec:fourier}.
The matrix $A$ has eigenvectors $(\chi_i(j))_{j\in\Gamma}$ for each
$i\in\Gamma$ with corresponding
eigenvalue $\wh{\beta}(i)$.
 Define $\wh{A}$ to be the matrix with $(i,j)$-entry
the Fourier transform $\wh{\beta}(i-j)$. For example, $\wh{I}=J$ and $\wh{J}=nI$.

Let us call a matrix {\em $\Gamma$-circulant} if, like the matrix $A$, it takes the form $(\alpha(i-j))$ for some function $\alpha\in\mathbb{C}^\Gamma$.
(When $\Gamma$ is cyclic these are circulant matrices as ordinarily defined.)
\begin{lem}\label{lem:isomorphism algebras}
The map $A\mapsto\wh{A}$ is an algebra isomorphism of the subalgebra
of $\Gamma$-circulant matrices under addition and ordinary matrix
multiplication with the subalgebra of $\Gamma$-circulant matrices
under addition and entrywise product:
$$\wh{A_1A_2}=\wh{A}_1\circ \wh{A}_2,\;\; \wh{A_1\circ
  A_2}=\frac{1}{n}\wh{A}_1\wh{A}_2.$$
In particular
$$\wh{A^k}=\wh{A}^{\;\circ k}, \;\; \wh{A^{\circ k}}=\frac{1}{n^{k-1}}\wh{A}^{\;k}.$$
\end{lem}
\begin{pf}
This is a consequence of the identities~\eqref{convolution} and \eqref{pointwise} for the Fourier transform of convolutions and pointwise products.
Since the Fourier transform is linear it follows that $A\rightarrow\wh{A}$ is an algebra isomorphism.
\end{pf}


\begin{lem}\label{lem:polynomials}
Given $$q(A)=q_0I+q_1A+q_1A^2+\cdots + q_\ell A^\ell,$$
define 
$$q^{\circ}(A)=q_0J+q_1A+q_2A^{\circ 2}+\cdots +q_\ell A^{\circ
 \ell},$$
i.e., the same polynomial except in Hadamard product powers of $A$ instead of
 ordinary matrix powers of $A$ and with the identity $I$ for ordinary
 matrix multiplication replaced by the identity $J$ for entrywise multiplication.

Then $\wh{q(A)}=q^{\circ}(\wh{A})$ and $\wh{q^{\circ}(A)}=q(\frac{1}{n}\wh{A})$.
\end{lem}
\begin{pf}
 This follows from Lemma~\ref{lem:isomorphism algebras}.
\end{pf}

Recall that for $Z\in\mathcal{G}_2$ to be a contractor for ${\rm hom}(\,\cdot\,, H)$ we need $M(Z)=I$.
The following is the key to obtaining such a contractor when $H$ has a $\Gamma$-circulant adjacency matrix $A$.
\begin{lem}\label{cor:id}
Suppose $A$ has eigenvector $\mathbf{1}$ with eigenvalue $s$, minimum
polynomial $p_A(t)$ and that $q(t)=p_A(t)/(t-s)$.
Then
$$I=\frac{1}{q(s)}q^{\circ}(\wh{A}).$$
\end{lem}
\begin{pf} By Corollary \ref{cor:J in powers of A} we have $J=\frac{n}{q(s)}q(A)$ and by  Lemma~\ref{lem:polynomials},  $\wh{q(A)}=q^{\circ}(\wh{A})$.
 Since $\wh{J}=nI$, the result follows. 
\end{pf}
Our main theorem of this section now results:
\begin{thm}\label{thm:contractor from Fourier transform}
Suppose $H$ is a connected graph with adjacency matrix a $\Gamma$-circulant matrix $A$ and that
$\wh{H}$ has adjacency matrix $\wh{A}$. Suppose further that $A$ has eigenvector
$\mathbf{1}$ with eigenvalue $s$, has minimum polynomial
$p_A(t)$, and that $$p_A(t)/(t-s)=q(t)=q_0+q_1t+\cdots + q_{\ell-1} t^{\ell-1}.$$

Then a contractor for ${\rm hom}(\,\cdot\,, \wh{H})$ is given by 
$$Z=\frac{1}{q(s)}\sum_{0\leq k\leq \ell-1} q_k P_1^{k},$$
where $P_1^{k}=C_k^*$ comprises two  vertices, one labelled $1$ and the other
$2$, joined by $k$ edges. (In particular, $P_1^{0}=\ol{K}_2$.)
\end{thm}
\begin{pf}
This is immediate from Lemma~\ref{cor:id} and the fact that $M_{\wh{H}}(P_1^{k})=\wh{A}^{\circ k}$. [The matrix $\wh{A}$ is the adjacency matrix of the edge-weighted graph $\wh{H}$, and $M_{\wh{H}}(X)$ the $V(\wh{H})\times V(\wh{H})$ matrix whose $(i,j)$-entry is ${\rm hom}_\phi(X,\wh{H})$, where $\phi(1)=i, \phi(2)=j$.]
\end{pf}

\section{A contractor for the number of $B$-flows}\label{sec:contractor_Sflows}

In this section we apply Theorem~\ref{thm:contractor from Fourier transform} to the special case where
 $H={\rm Cayley}(\Gamma, B)$ for additive Abelian group $\Gamma$, and where $B\subseteq \Gamma\setminus 0$ satisfies $-B=B$. 

Let
$\delta_B$ be the indicator function of $B$. The
adjacency matrix $A$ of $H$ has $(i,j)$-entry $\delta_B(i-j)$ and eigenvalues $\wh{\delta}_B(i)$ for $i\in\Gamma$.
The characteristic polynomial of $A$ is given by 
$$ \phi_A(t)=\prod_{c\in\Gamma}(t-\wh{\delta_B}(c)),$$
and has degree $|\Gamma|=|V(H)|$. The minimum polynomial of $A$ has
degree $\ell$ equal to the number of distinct values of $\wh{\delta_B}(c)$
for $c\in \Gamma$, bounded below by $d(H)+1$, where $d(H)$ is the
diameter of $H$. (This lower bound is attained when $H$ is a distance
regular graph.)
The largest eigenvalue of $A$ is $\wh{\delta_B}(0)=|B|$, belonging to the
eigenvector $\mathbf{1}$.

We reach the main theorem of this paper, giving a contractor for $B$-flows.
\begin{thm}\label{thm.S-flows}
Let $\Gamma$ be an additive Abelian group of order $n$ and suppose $H={\rm Cayley}(\Gamma, B)$ is a connected graph with adjacency matrix $A$, whose
$(i,j)$ entry is equal to $\delta_B(i-j)$. Let $A$ have minimum
polynomial 
$p_A(t)$
 of degree $\ell$ equal to the number of distinct
eigenvalues $\wh{\delta_B}(i)$ of $A$, and set 
$$p_A(t)/(t-|B|)=q(t)=q_0+q_1t+\cdots + q_{\ell-1}
t^{\ell-1}.$$ 
Define $\wh{H}$ to be the
edge-weighted graph on vertex set $\Gamma$ with adjacency matrix
having $(i,j)$-entry $\wh{\delta_B}(i-j)$.
Then

\begin{itemize}


\item[(i)] A contractor for the number of $B$-flows is given by 
$$Z=\frac{n}{q(|B|)}\sum_{0\leq k\leq \ell-1} q_k P_1^{k},$$
where $P_1^{k}=C_k^*$ comprises two  vertices, one labelled $1$ and the other
$2$, joined by $k$ parallel edges (in particular, $P_1^0=\ol{K}_2$);

\item[(ii)] A contractor for $n^{c(G)}$ times the number of $B$-tensions of $G$
is given by
$$Z=\frac{1}{q_0}\left[\frac{q(|B|)}{n}\ol{K}_2-\sum_{1\leq k\leq \ell-1}q_kP_{k}\right].$$
\end{itemize}
\end{thm}

\begin{pf}
Part (i) is a corollary of Theorem~\ref{thm:contractor from Fourier transform}, which gives 
$$\frac{1}{q(|B|)}\sum_{0\leq k\leq \ell-1} q_k P_1^{k},$$ 
as a contractor for ${\rm hom}(\,\cdot\, , \wh{H})$. The number of $B$-flows of $G$ is by Lemma~\ref{lem:Sflows} equal to $n^{-|V(G)|}{\rm hom}(G,\wh{H})$, and this implies the expression for the contractor as given in the theorem statement.

In order to show (ii), by Corollary~\ref{cor:J in powers of A} we have
$$J=\frac{n}{q(|B|)}\sum_{0\leq k\leq \ell-1}q_kA^k.$$
Since $q_0\neq 0$, we then obtain
$$I=\frac{1}{q_0}\left[\frac{q(|B|)}{n}J-\sum_{1\leq k\leq \ell-1}q_kA^k\right],$$
and this implies $M_H(Z)=I$, where $Z$ is the quantum graph given in the theorem statement, i.e., $Z$ is a contractor for ${\rm hom}(\, . \,, H)$. By Lemma~\ref{lem:Sflows}, the number of $B$-tensions of $G$ is equal to $n^{-c(G)}{\rm hom}(G,H)$.
\end{pf}

\section{Examples}\label{sec:ex}

In this section we exhibit some concrete instances of the contractor for $B$-flows given by Theorem~\ref{thm.S-flows}.

\subsection{Proper vertex colourings, nowhere-zero flows}

When $\Gamma$ is order $n$ and $B=\Gamma\setminus\{0\}$, the graph $H={\rm Cayley}(\Gamma, B)$ is
isomorphic to $K_n$. The adjacency matrix $A=J-I$ has minimum
polynomial $p_A(t)=(t+1)(t-n+1)$ and $q(t)=p_A(t)/(t-n+1)=1+t$. Here $J=I+A$.
A contractor for ${\rm hom}(\,\cdot\, , K_n)$ is given by $\ol{K}_2-P_1$
(which amounts to the contraction--deletion identity for the chromatic
polynomial).
Since $A=(n\!-\!1)J-A^2$, a connector for ${\rm hom}(\,\cdot\, , K_n)$ is
given by $(n-1)\ol{K}_2-P_2$. (In~{\rm \cite{LS09}} the connector given is
$\frac{1}{n-1}\left[P_3-(n-2)P_2\right]$, which follows from the
equation $(n-1)A=A^3-(n-2)A^2$ obtained from $Ap_A(A)=0$.) 

The graph $\wh{H}$ has adjacency matrix $\wh{A}=nI-J$, which satisfies
$\wh{A}^2=n\wh{A}$, so that $\frac{1}{n}P_2$ is a connector for ${\rm hom}(G,
\wh{H})=n^{|V(G)|}F(G;n)$, and hence $P_2$ is a connector for $F(G;n)$, the number of nowhere-zero $\mathbb{Z}_n$ flows (as expected).
Since $I=\frac{1}{n}(J+\wh{A})$, a contractor is given by
$\frac{1}{n}\left[\,\ol{K}_2+P_1\right]$. Hence a contractor for $F(G;n)$ is given by $\ol{K}_2+P_1$. (That this is the case amounts to the deletion--contraction identity for the flow polynomial).

\subsection{Strongly regular graphs}
If $H$ is a strongly regular graph with parameters $(n,k,\lambda,\mu)$
then its adjacency matrix satisfies $$A^2+(\mu-\lambda)A+(\mu-k)I=\mu
J.$$
In this case, a contractor $Z$ for ${\rm hom}(\,\cdot\, ,
H)$ is thus given by
$$(\mu-k)Z=\mu\ol{K}_2-(\mu-\lambda)P_1-P_2.$$
This is the content of {\rm \cite[Proposition 8]{dlHJ95}}. 

The graph ${\rm Cayley}(\Gamma,B)$ is a strongly regular graph with parameters $(n,k,\lambda,\mu)$ precisely
when $B$ is a $(n,k,\lambda,\mu)$-partial difference set (the differences $i-j$
for $i,j\in B$ give each element of $B$ exactly $\lambda$ times, each element of
$\Gamma-B-\{0\}$ appears $\mu$ times, and $0$ appears $k$ times). For
example, the Paley graph of order $q\equiv 1\pmod 4$
($\Gamma=\mathbb{F}_q$, $B$ the set of non-zero squares) is a strongly regular graph with parameters $(q,\frac{q-1}{2},\frac{q-5}{4},\frac{q-1}{4})$.


\subsection{An intermediate result}

In order to proceed with examples related to some of the $B$-flow conjectures described in~\cite{Jae88}, we need to supplement Theorem~\ref{thm.S-flows} by its statement for disconnected Cayley graphs $H$.

If $B$ does not contain a set of generators for $\Gamma$ then the $\Gamma$-circulant graph $H={\rm Cayley}(\Gamma, B)$ is not connected. The connected components of ${\rm Cayley}(\Gamma, B)$ correspond to the cosets of the subgroup generated by $B$. 
\begin{thm}\label{cor:disconn}
Let $\Gamma$ be an additive Abelian group of order $n$ and suppose $H={\rm Cayley}(\Gamma, B)$ has adjacency matrix $A$. 
Let $A$ have minimum
polynomial 
$p_A(t)$
and set
 $$p_A(t)/(t-|B|)=q(t)=q_0+q_1t+\cdots + q_{\ell-1}
t^{\ell-1}.$$ 

Suppose $H$ has $a$ isomorphic connected components, i.e., its adjacency matrix $A$ is permutation-equivalent to a matrix of the form $I\otimes A_1$ for $a\times a$ identity matrix $I$ and some $n/a\times n/a$ matrix $A_1$.
Then, a contractor for the number of $B$-flows is given by 
$$Z=\frac{n}{aq(|B|)}\sum_{0\leq k\leq \ell-1} q_k P_1^{k},$$
where $P_1^{k}$ consists of two  vertices, one labelled $1$ and the other
$2$, joined by $k$ edges. 
\end{thm}
\begin{pf} If $B$ generates a subgroup $\Gamma_1$ of index $a$ in $\Gamma$ then ${\rm Cayley}(\Gamma,B)$ has $a$ components, each isomorphic to a connected Cayley graph $H_1={\rm Cayley}(\Gamma_1, B)$ of order $n_1=|\Gamma_1|=n/a$. Then ${\rm hom}(G,H)=a^{c(G)}{\rm hom}(G,H_1)$ is equal to the number of $B$-tensions of $G$ multiplied by $n^{c(G)}$, and $a^{c(G)}n_1^{-|V(G)|}{\rm hom}(G,\wh{H}_1)$ is equal to the number of $B$-flows of $G$. If $H_1$ has adjacency matrix $A_1$ then $H$ has adjacency matrix $I\otimes A_1$, where $I$ is the $a\times a$ identity matrix. Let $q_1(t)=p_{A_1}(t)/(t-|B|)$, where $p_{A_1}(t)$ is the minimum polynomial of $A_1$. Then $q_1(t)=q(t)$, since the eigenvalues of $A_1$ coincide with those of $A$ (only the eigenspaces differ). By Corollary~\ref{cor:J in powers of A}, $J=\frac{n_1}{q_1(s)}q_1(A_1)$, whence $I=\frac{n_1}{q(|B|)}q^\circ(\wh{A_1})$ by Corollary~\ref{cor:id}.  The result now follows in the same fashion as for the proof of Theorem~\ref{thm.S-flows}. 
\end{pf}

\subsection{Cycle double covers}
Let $\Delta$ be an additive group of $m$ elements, and $\Delta^d$ the $d$-fold Cartesian product. 
The Hamming weight $|x|$ of a given element $x=(x_1,\ldots, x_d)\in \Delta^d$ is defined by $|x|=\#\{i:x_i\neq 0\}$.
The sets $S_{r}^{d}=\{x\in \Delta^d:|x|=r\}$ are called {\em shells} (of radius  $r$).
 The {\em Krawtchouk polynomial} $K_r(w;d,m)$ of degree $r$ is defined  for $0\leq r, w\leq d$ by 
$$K_r(w;d,m)=[z^r](1+(g-1)z)^{d-w}(1-z)^w=\sum_{0\leq i\leq r}(-1)^i(m-1)^{r-i}\binom{w}{i}\binom{d-w}{r-i}.$$
In particular, $K_0(w;d,m)=1, K_1(w;d,m)=(m-1)d-mw$ and $K_d(w;d,m)=(-1)^{w}(m-1)^{d-w}$.
The Fourier transform of the indicator function $\delta_{S_{r}^d}$ of a shell is given (see e.g.~\cite{HP03}) by 
$$\wh{\delta_{S_{r}^d}}(x)=K_r(|x|;d,m).$$ 

When $m=2$ ($\Delta=\mathbb{F}_2$),
$$K_r(w;d,2)=\sum_{0\leq i\leq r}(-1)^i\binom{w}{i}\binom{d-w}{r-i}.$$

Fix $m=2$ for now.

We begin with a trivial illustration of Theorem~\ref{cor:disconn}.

Let $S_d^d\subseteq\mathbb{F}_2^d$ be the single-element subset comprising the all-one vector of Hamming weight $d$. Then ${\rm Cayley}(\mathbb{F}_2^d,S_2^d)$ consists of $2^{d-1}$ copies of $K_2$. Here $\wh{\delta_{S_d^d}}(x)=K_d(|x|;d,2)=(-1)^{|x|}$, from which we find $q(t)=t+1$, and a contractor for $S_d^d$-flows is given by Theorem~\ref{cor:disconn} as 
$$\frac{2^d}{2^{d-1}q(1)}\left[P_1^1+P_1^0\right]=P_1+\ol{K}_2.$$
This, as it should be, is also a contractor for the number of $\mathbb{F}_2$-flows.

Let $F_{S_{r}^d}(G)$ denote the number of $S_r^d$-flows of $G$, where $S_{r}^d$ is the shell of radius $r$ in $\mathbb{F}_2^d$, i.e., the set of vectors in $\mathbb{F}_2^d$ with exactly $r$ non-zero coordinates.
 
In~\cite{Jae88} the Double Cover Conjecture is reformulated as a $B$-flow problem:
\begin{problem} Is it true that for every bridgeless graph there exists $d\geq 2$ so that $F_{S_2^d}(G)\neq 0$?
\end{problem}
Jaeger also mentions the stronger 5-colourable Cycle Double Cover Conjecture, that in fact $d=5$ suffices, i.e., that $F_{S_2^5}(G)\neq 0$ for every bridgeless graph $G$.

We have just seen that a graph $G$ has $F_{S_2^2}(G)\neq 0$ if and only if $G$ has a nowhere-zero $\mathbb{F}_2$-flow. 
It is a result due originally to Tutte that $F_{S_2^3}(G)\neq 0$ if and only if $G$ has a nowhere-zero $\mathbb{F}_4$-flow. 

\subsubsection{3-colourable cycle double covers}

Let $S_2^3\subseteq \mathbb{F}_2^3$. Then $H={\rm Cayley}(\mathbb{F}_2^3,S_2^3)$ consists of two connected components isomorphic to $K_4$. (Hence there are $S_{2}^3$-tensions of $G$ if and only if there are $\mathbb{F}_4$-tensions of $G$.) The adjacency matrix of $H$ has eigenvalues $\wh{\delta}_{S_2^3}(x)\in\{3,-1\}$, giving $q(t)=t+1$. 
A contractor for $S_2^3$-flows is by Theorem~\ref{cor:disconn} given by 
$$\frac{8}{2\cdot q(3)}[P_1^1+P_1^0]=P_1+\ol{K}_2.$$
Again, this is as it should be as $S_2^3$-flows are in one-to-one correspondence with nowhere-zero $\mathbb{F}_4$-flows (truncate the elements of $S_2^3\subseteq\mathbb{F}_2^3$ in the last position; conversely, to a nowhere-zero $\mathbb{F}_2^2$-flow add a parity check digit to make a $S_2^3$-flow).

\subsubsection{4-colourable cycle double covers}


 Let $S_2^4\subseteq \mathbb{F}_2^4$. Then $H={\rm Cayley}(\mathbb{F}_2^4,S_2^4)$ consists of two isomorphic connected components, each a regular graph of degree $\binom{4}{2}=6$ on $8$ vertices (i.e., $K_8$ minus a perfect matching). The adjacency matrix of $H$ has eigenvalues $\wh{\delta}_{S_2^4}(x)\in\{6,0,-2\}$, giving $q(t)=t(t+2)$. 
A contractor for $S_2^4$-flows is by Theorem~\ref{cor:disconn} given by 
$$\frac{16}{2\cdot q(6)}[P_1^2+2P_1^1]=\frac{1}{6}[C_2+2P_1].$$

\subsubsection{$d$-colourable cycle double covers}

We now consider the general case of $d$-colourable cycle double covers. 
In this case we have
\begin{align*}\wh{\delta_{S_2^d}}(x) & =K_2(|x|;d,2)\\
 & = \binom{d-|x|}{2}-|x|(d-|x|)+\binom{|x|}{2}\\
 & = \binom{d}{2}+2|x|(|x|-d).\end{align*}
The adjacency matrix $A$ of the Cayley graph $H={\rm Cayley}(\mathbb{F}_2^d, S_{2}^d)$ has $\lceil \frac{d+1}{2}\rceil$ distinct eigenvalues, namely $\binom{d}{2}+2w(w-d)$ for $0\leq w\leq \lfloor \frac{d}{2}\rfloor$. The largest eigenvalue $\binom{d}{2}=|S_{2}^d|$ belongs to the eigenvector $\mathbf{1}$.
In the notation of Theorem~\ref{cor:disconn},
$$q(t)=\prod_{1\leq w\leq \lfloor d/2\rfloor}\left[t-\binom{d}{2}+2w(d-w)\right],$$
and 
$$q(\binom{d}{2})=\begin{cases}2^{d/2-1}d! & \mbox{\rm $d$ even,}\\
 2^{(d-1)/2}(d-1)! & \mbox{\rm $d$ odd.}\end{cases}$$
Since the elements of $S_2^d$ generate precisely the even-weight vectors in $\mathbb{F}_2^d$, the graph $H={\rm Cayley}(\mathbb{F}_2^d,S_2^d)$ has two connected components.
This allows us to explicitly construct a contractor for $F_{S_2^d}(G)$, as we have done for $d\in\{2,3,4\}$ and we now do for the case $d=5$.

\subsubsection{$5$-colourable cycle double covers}

Let $A$ be the adjacency matrix of ${\rm Cayley}(\mathbb{F}_2^5, S_2^5)$, which has $a=2$ connected components, eigenvalues $10, 2,-2$, and $q(t)=(t-2)(t+2)=t^2-4$, i.e., $q_0=-4, q_1=0, q_2=1$, and $q(10)=96$.
We have $n=2^5=32$ and $F_{S_2^5}(G)=32^{-|V(G)|}{\rm hom}(G,\wh{H})$. Hence, by Theorem~\ref{cor:disconn} a contractor for $F_{S_2^5}(G)$ is given by
$$\frac{32}{2\cdot q(10)}[P_1^2-4P_1^0]=\frac{1}{6}\left[C_2-4\ol{K}_2\right].$$
If we let $e$ be an edge with endpoints labelled $1$ and $2$ and we write $G=XP_1$, then $X=G\backslash e$ and $XK_1=G/e$. Suppose we write $G\| e$ for the operation that inserts an extra edge parallel to $e\in E(G)$. Then $XP_1^2=G\| e$.
In this notation,
$$F_{S_2^5}(G\| e)=4F_{S_2^5}(G\backslash e)+6F_{S_2^5}(G/e).$$


\subsection{Fulkerson flows}

Fulkerson's Conjecture is that in every bridgeless cubic graph $G$ there is a family of six perfect matchings such that each edge appears in exactly two of them.
Jaeger~\cite[Theorem 6.1]{Jae88} shows that Fulkerson's Conjecture is equivalent to the assertion that a bridgeless cubic graph has a $S_4^6$-flow, where $S_4^6\subseteq\mathbb{F}_2^6$ comprises those vectors containing exactly four 1s.
We follow~\cite{DVNR02} and call these types of flow ``Fulkerson flows''.
 
Let $A$ be the adjacency matrix of ${\rm Cayley}(\mathbb{F}_2^6, S_4^6)$, which has $a=2$ connected components, eigenvalues $15, -5,-1,3$, and $q(t)=(t+5)(t+1)(t-3)=t^3+3t^2-13t-15$, i.e., $q_0=-15, q_1=-13, q_2=3$, $q_3=1$ and $q(15)=20\cdot 16\cdot 12$, $n=2^6$.
We have $F_{S_4^6}(G)=64^{-|V(G)|}{\rm hom}(G,\wh{H})$. Hence, by Theorem~\ref{cor:disconn} a contractor for $F_{S_4^6}(G)$ is given by
$$\frac{64}{2\cdot q(15)}[P_1^3+3P_1^2-13P_1^1-15P_1^0]=\frac{1}{120}\left[P_1^3+3C_2-13P_1-15\ol{K}_2\right],$$
where $P_1^3=C_3^*$ consists of two vertices joined by 3 parallel edges.

After similar calculations, a contractor for $F_{S_2^6}(G)$ (the next case in the Double Cover Conjecture series of examples begun above) is given by
$$\frac{1}{90}\left[P_1^3-C_2-17P_1-15\ol{K}_2\right].$$

\subsection{Petersen flows}\label{ex:petersen}

The Petersen Flow Conjecture of Jaeger~\cite{Jae88} states that every bridgeless graph has a Petersen flow. As explained in e.g.~\cite{DVNR02}, this conjecture implies the 5-colourable cycle double conjecture and other long-standing conjectures.
Petersen flows are defined as follows. 

Take $\Gamma=\mathbb{F}_2^6$.
Let $C^{(1)}, \ldots , C^{(6)}$ be a basis of the cycle space of the Petersen graph $P=(V,E)$ and let $y^{(1)}, \ldots, y^{(6)}$ be the corresponding indicator vectors in $\mathbb{F}_2^{E}$ (where $|E|=15$). Let $B\subseteq\mathbb{F}_2^6$ be the set of vectors $(y_e^{(1)}, \ldots, y_e^{(6)})$ for each $e\in E$. Then $|B|=15$ since for every pair of distinct edges there is a cycle in the basis which contains exactly one of them.

The indicator function $\delta_B$ has Fourier transform
\begin{align*}\wh{\delta_B}(a_1\ldots a_6) &=\sum_{b_1\ldots b_6\in B}(-1)^{a_1b_1+\cdots + a_6b_6}\\
& = \#\{e: (a_1y^{(1)}+\cdots + a_6y^{(6)})_e=0\} -\#\{e: (a_1y^{(1)}+\cdots + a_6y^{(6)})_e=1\}\\
& = |E|-2|C|,\end{align*}
where $C$ is the cycle with indicator vector $a_1y^{(1)}+\cdots + a_6y^{(6)}$.
Since the Petersen graph has cycles of lengths $0, 5, 6, 8$ and $9$, it follows that ${\rm Cayley}(\mathbb{F}_2^6, B)$ has eigenvalues $15, 5,3,-1$ and $-3$.
Thus 
$$q(t)=(t+1)(t^2-9)(t-5)=t^4-4t^3-14t^2+36t+45,$$
and $n=64, q(15)=34560$. Since $B$ contains a set of generators for $\mathbb{F}_2^6$ the graph ${\rm Cayley}(\mathbb{F}_2^6, B)$ is connected. 
Hence a contractor for the number of $B$-flows is given by
$$\frac{64}{34560}\left[P_1^4-4P_1^3-14P_1^2+36P_1+45\ol{K}_2\right],$$
the fraction simplifying to $\frac{1}{540}$.

This construction of a contractor for Petersen flows extends more generally to $B$-flows when $B$ is defined in a similar way via a basis for the cycle space of a graph $Q=(V,E)$, the eigenvalues of ${\rm Cayley}(\mathbb{F}_2^E,B)$ equal to $|E|-2|C|$ for Eulerian subgraphs $C$ of $Q$. ($B$-flows of $G$ here are cycle-continuous maps from $G$ to $Q$. See e.g.~\cite{DVNR02}.)   

\subsection{Summary}

\begin{figure}[htp]
\centering
\caption{Examples of contractors for a selection of graph invariants. The identities are as seen modulo the graph parameter in question. The first example is to be understood as saying that $\ol{K}_2-K_2$ is a contractor for the chromatic polynomial, i.e., $P(XK_1;q)=P(X\ol{K}_2;q)-P(XK_2;q)$ for every graph $X$ with a vertex labelled $1$ and a vertex labelled $2$, where $K_1$ has labels $1$ and $2$ on its single vertex, $K_2$ has one vertex label $1$ and the other label $2$, and similarly for $\ol{K}_2$. The graph $XK_1$ is the result of identifying the vertices labelled $1$ and $2$ in $X$, $XK_2$ adds an edge between $1$ and $2$, etc. 
}\label{fig:contractors}
\bigskip
\hspace{-0.4cm}\includegraphics[scale=0.45]{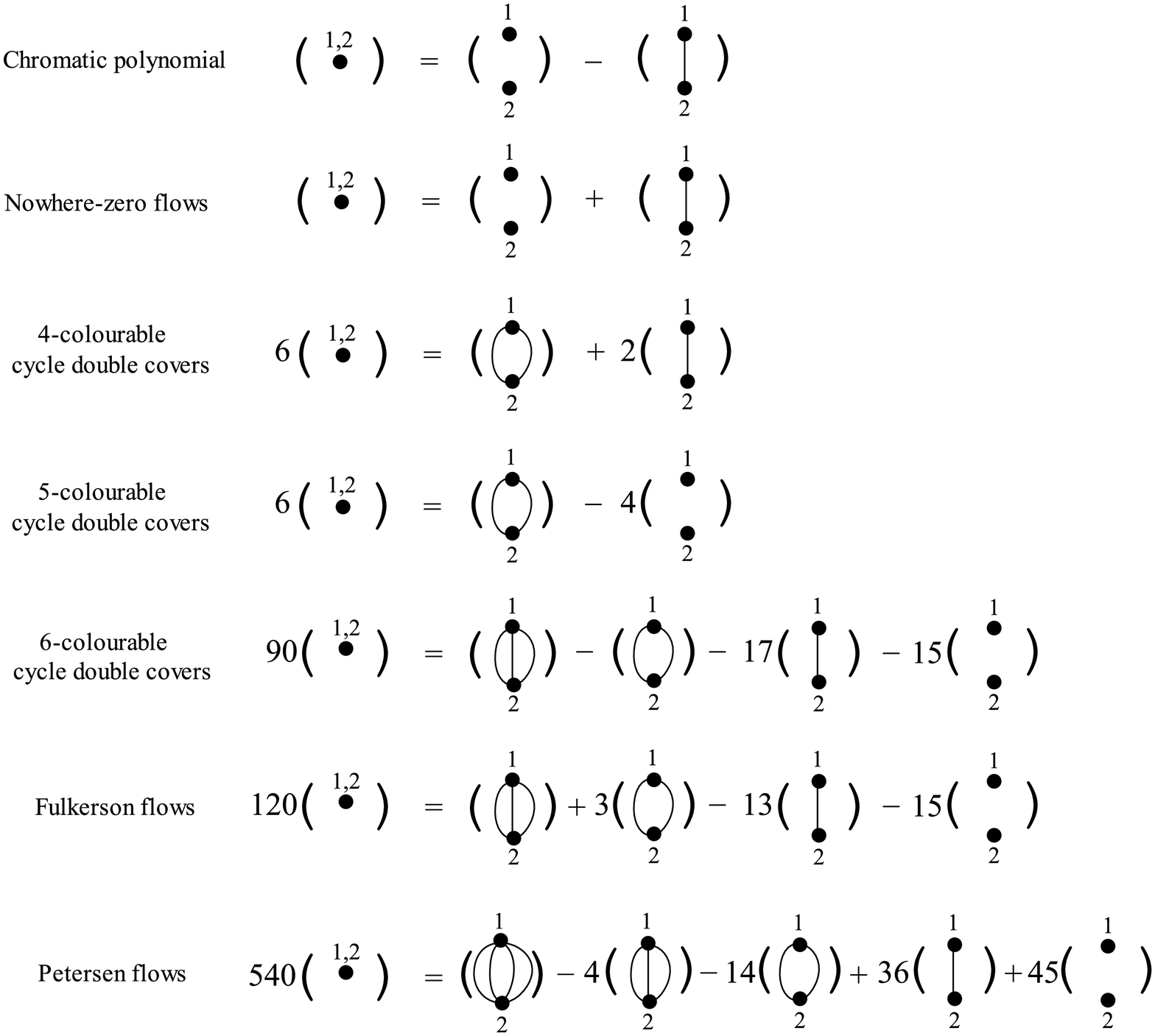}

\end{figure}

In Figure~\ref{fig:contractors} we display contractors for those $B$-flows that we have encountered so far.

In general, to find a contractor for the number of $B$-flows we need to find the $\ell$ distinct eigenvalues of the adjacency matrix of ${\rm Cayley}(\Gamma, B)$ ($\wh{\delta_{B}}(x)$ for $x\in\Gamma$), discard the largest one ($|B|$), and then the $k$th elementary symmetric functions of the remaining $\ell-1$ eigenvalues form the coefficients in the expression for a contractor as a linear combination of graphs $P_1^{\ell-1-k}$ (two vertices joined by $\ell-1-k$ edges).

We finish our selection of examples with one that is of interest because of its connection with the graph parameter counting Eulerian orientations.

\subsection{Nowhere-zero $n,1$-flows}\label{ex:Eulerian}  

Take $\Gamma=\mathbb{Z}_n$, the integers modulo $n$, and $B=\{-1,+1\}$. A $B$-flow in a graph $G$ is an orientation of $G$ in which the indegree is congruent to the outdegree modulo $n$ (also known as a nowhere-zero $n,1-$flow).   
When $n=3$ we have $B$-flows equal to nowhere-zero flows, and a contractor is given by $\ol{K}_2+K_2$. In case $G$ is a 4-regular graph,  $B$-flows of $G$ are Eulerian orientations of $G$. 
More generally, if we take $B$-flows modulo $n$ and all the vertex degrees of $G$ belong to $\{0,1,2,\ldots, n-1, n+1\}$ then $B$-flows modulo $n$ correspond to Eulerian orientations of $G$.

The indicator function of $B$ has Fourier transform $\wh{\delta_B}(a)=\zeta^a+\zeta^{-a}$, where $\zeta=e^{2\pi i/n}$.
When $n=4$ these eigenvalues are $2, 0, -1$ and we find that $q(t)=t(t+2)$, so that $P_1^2+2P_1$ is a contractor for $B$-flows modulo $4$ (where we recall that $P_1^k$ denotes the graph with two vertices joined by $k$ edges).
When $n=5$ we have 
$$q(t)=(t-\zeta-\zeta^{-1})(t-\zeta^2-\zeta^{-2})=t^2+t-1,$$
so $P_2+P_1-\ol{K}_2$ is a contractor for $B$-flows modulo $5$.

Generally then, for $B$-flows modulo $n$ we have
$$q(t)=\prod_{a=1}^{\lfloor n/2\rfloor}(t-\zeta^a-\zeta^{-a}).$$
When $n$ is odd, 
$$q(2)=\prod_{a=1}^{\frac{n-1}{2}}(1-\zeta^{a})(1-\zeta^{-a})=n,$$
and when $n$ is even $q(2)=2n$.

Take $n$ odd and let $(-1)^{j}e_j$ denote the coefficient of $t^{\frac{n-1}{2}-j}$ in $q(t)$. Let 
$$s_k=\sum_{a=1}^{\frac{n-1}{2}}(\zeta^a+\zeta^{-a})^k,$$
denote the $k$th power sum symmetric function in the roots of $q(t)$.
By the Newton-Girard identities the elementary symemtric functions in the roots of $q(t)$ are given by
$$e_j= \sum_{k_1+2k_2+\cdots + jk_j=j}\frac{(-1)^{k_1+k_2+\cdots + k_j}}{k_1!k_2!\cdots k_j!}\left(\frac{s_1}{1}\right)^{k_1}\left(\frac{s_2}{2}\right)^{k_2}\cdots\left(\frac{s_j}{j}\right)^{k_j},$$
the sum over all partitions of $j$ with $k_i$ parts equal to $i$.
We find that
\begin{align*}s_k & = \sum_{b=1}^{n-1}\sum_{j=0}^{\lfloor \frac{k}{2}\rfloor}\binom{k}{j}\zeta^{(k-2j)b}\\
& \begin{cases} = -\sum_{j=0}^{\frac{k-1}{2}}\binom{k}{j} = -2^{k-1} & \mbox{$k$ odd,}\\
= -\sum_{j=0}^{\frac{k}{2}-1}\binom{k}{j}+\frac{n-1}{2}\binom{k}{k/2}= \frac{n}{2}\binom{k}{k/2}-2^{k-1} & \mbox{$k$ even.}\end{cases}\end{align*}
So $s_1=-1, s_2=n-2, s_3=-4, s_4=3n-8, s_5=-16$ etc., giving
$e_1=-1$, $e_2=-\frac{n-3}{2}$ for odd $n\geq 3$, $e_3=\frac{n-5}{2}$ for odd $n\geq 5$, $e_4=\frac{(n-5)(n-7)}{8}$ for odd $n\geq 5$.
Thus for odd $n$ the polynomial $q(t)$ begins
$$q(t)=t^{\frac{n-1}{2}}+t^{\frac{n-1}{2}-1}-\frac{n-3}{2}t^{\frac{n-1}{2}-2}-\frac{n-5}{2}t^{\frac{n-1}{2}-3}+\frac{(n-5)(n-7)}{8}t^{\frac{n-1}{2}-4}+\cdots$$
and the contractor for $B$-flows modulo odd $n$ accordingly begins
$$P_1^{\frac{n-1}{2}}+P_1^{\frac{n-1}{2}-1}-\frac{n-3}{2}P_1^{\frac{n-1}{2}-2}-\frac{n-5}{2}P_1^{\frac{n-1}{2}-3}+\frac{(n-5)(n-7)}{8}P_1^{\frac{n-1}{2}-4}+\cdots$$
Is there a simple formula for each of the coefficients of $q(t)$?

In Table~\ref{table:Eulerian_orns} we list contractors for $\{\pm 1\}$-flows modulo $n$ for $3\leq n\leq 9$. 
For graphs whose vertex degrees all belong to $\{0,1,2,\ldots, n-1, n+1\}$, the number of $\{\pm 1\}$-flows is equal to the number of Eulerian orientations of $G$.
\begin{table}[ht]\caption{}\label{table:Eulerian_orns}
\vspace{0.2cm}
\centering
\begin{tabular}{l|l|l}
\hline\hline
$n$ & $q(t)$ & Contractor for $\{\pm 1\}$-flows modulo $n$  \\[0.5ex]
\hline
$3$ & $1+t$ & $\ol{K}_2+P_1$ \\
$4$  & $2t+t^2$ &$P_1+\frac12P_1^2$ \\
$5$ & $-1+t+t^2$ & $-\ol{K}_2+P_1+P_1^2$\\
$6$ & $-2-t+2t^2+t^3$ & $-\ol{K}_2-\frac12P_1+P_1^2+\frac12P_1^3$ \\
$7$ & $-1-2t+t^2+t^3$ & $-\ol{K}_2-2P_1+P_1^2+P_1^3$ \\
$8$ & $-4t-2t^2+2t^3+t^4$ & $-2P_1-P_1^2+P_1^3+\frac12 P_1^4$\\
$9$ & $1-2t-3t^2+t^3+t^4$ & $\ol{K}_2-2P_1-3P_1^2+P_1^3+P_1^4$\\[1ex]
\hline
\end{tabular}

\end{table} 
 
We have seen that the number of Eulerian orientations is expressible as the limit of the graph parameters counting the number of $\{\pm 1\}$-flows in $\mathbb{Z}_n$ as $n\rightarrow\infty$. Equivalently, Eulerian orientations are counted by the number of homomorphisms to the Fourier dual of the infinite graph ${\rm Cayley}(\mathbb{Z},\{\pm 1\})$, which is the graph on vertex set the real interval $[0,1]$ with edges $ij$ weighted $2\cos 2\pi (i-j)$, equal to the Fourier transform of the edge weights on ${\rm Cayley}(\mathbb{Z}, \{\pm 1\})$. (The dual group of $\mathbb{Z}$ is the circle group, which in the usual way can be identifed with the additive group of reals in $[0,1]$ modulo $1$.) This is the limiting graph density for the graph parameter counting Eulerian orientations that is described in~\cite{LS06, LS09}. (In the notation of these two papers, in terms of ``graphons'' the number of Eulerian orientations of $G$ can be expressed as $t(G,W)$, where $W(x,y)=2\cos(2\pi(x-y))$.)  
Is there any sense in which we can say these contractors for $\{\pm 1\}$-flows modulo $n$ converge to a ``limit contractor'' for Eulerian orientations?

\section{Conclusion}\label{sec:conclusion}

In this paper we provide a new and constructive proof for the existence of a contractor for $B$-flows. (Lov\'asz and B.~Szegedy~\cite{LS09} do not seem to provide an explicit construction in their proof of the existence of a series-parallel contractor for homomorphism numbers.) 

Theorem~\ref{thm:contractor from Fourier transform} provides  a construction of a contractor for the number of homomorphisms to an edge-weighted Cayley graph. Specifically, it gives a contractor for ${\rm hom}(\, . \, , \wh{H})$ when $H$ is an edge-weighted Cayley graph of a finite Abelian group $\Gamma$, and $\wh{H}$ is the edge-weighted Cayley graph with weights given by the Fourier transform of the weights in $H$; Theorem~\ref{thm.S-flows} applies this to the special case of $B$-flows, where $H={\rm Cayley}(\Gamma, B)$.

Let us finish this paper with a few problems.

Left open by Theorem~\ref{thm:connector contractor}, and only partially covered by Theorem~\ref{thm.S-flows}, is the following:
\begin{problem}\label{prob:contractor_eigenvalue_zero}
Given an edge-weighted connected graph $H$ with adjacency matrix $A$ that has $0$ as an eigenvalue, find a contractor for ${\rm hom}(\,\cdot\, , H)$.  
 \end{problem}

By the result of Lov\'asz and B.~Szegedy, in Problem  \ref{prob:contractor_eigenvalue_zero} it suffices to restrict the search to series-parallel graphs.

Recall that the quantum graph $Z$ is a contractor for the graph parameter ${\rm hom}(\,\cdot\, , H)$ if ${\rm hom}_\phi(G(Z-K_1),H)=0$ for each $G\in\mathcal{G}_2$ and $\phi:[2]\rightarrow V(G)$. This equality has (with a switch of homomorphisms from the left and from the right) a similar form to the inequality in the following problem: 
\be\label{eq:undecidable} \begin{array}{c} \mbox{Given a $k$-labelled quantum graph $Y$, is it the case that} \\ \mbox{\rm  ${\rm hom}_\phi(Y,G)\geq 0$ for all $k$-labelled graphs $G$ and maps $\phi: [k]\rightarrow V(G)$?}\end{array}\ee

Lov\'asz's Seventeenth Problem asks whether such a positive quantum graph $Y$ can be expressed as the sum of squares of labelled quantum graphs (with multiplication being the disjoint union of labelled graphs followed by identification of like-labelled vertices). This has been recently answered in the negative by Hatami and Norine~\cite{HN10}: Problem~\eqref{eq:undecidable} is in general undecidable. 
It would be interesting to know what the complexity of this problem is when $Y$ is restricted to a subclass of quantum graphs, such as series-parallel graphs, and likewise the computational complexity of Problem~\ref{prob:contractor_eigenvalue_zero}.  


Finally, we highlight a particular case of Problem~\ref{prob:contractor_eigenvalue_zero}, since it may be of interest to compare its solution to the three-term deletion--contraction--elimination recursion for the polynomial of Averbouch--Godlin--Makowsky~\cite{AGM07}. 
\begin{problem} 
  Find a contractor for the graph parameter ${\rm hom}(G,
  K_{q-p}^1+K_p^y)$, where $K_p^y$ denotes the complete graph on $p$ vertices with a loop of weight $y$ on each vertex (and similarly for $K_{q-p}^1$), and the $+$ denotes graph join. (This homomorphism number is equal to $\xi_G(q,y\!-\!1,(p\!-\!q)(y\!-\!1))$ in the notation of~{\rm \cite{AGM07}}). 
\end{problem}



\end{document}